\documentclass[12pt,twoside]{amsart}
\usepackage{amssymb}
\usepackage{amscd}

\title{$\mathbf{\overline{C}_{n,n-1}}$ revisited}
\author{Osamu Fujino} 
\date{2005/7/29}
\address{Graduate School of Mathematics\\ 
 Nagoya University, Chikusa-ku Nagoya 464-8602 Japan}
\email{fujino@math.nagoya-u.ac.jp}

\newcommand{\xhor}[0]{{\operatorname{hor}}}
\newcommand{\xSupp}[0]{{\operatorname{Supp}}}
\newcommand{\Exc}[0]{{\operatorname{Exc}}}
\newcommand{\Hom}[0]{{\operatorname{Hom}}}
\newcommand{\Var}[0]{{\operatorname{Var}}}
\newcommand{\red}[0]{{\operatorname{red}}}
\newcommand{\xver}[0]{{\operatorname{ver}}}
\newtheorem{thm}{Theorem}[section]
\newtheorem{lem}[thm]{Lemma}

\theoremstyle{definition}

\newtheorem*{ack}{Acknowledgments}       
 
\newtheorem*{notation}{Notation}         
       
\newtheorem*{case}{Case}         

\newtheorem{step}{Step}

\begin{document}
\bibliographystyle{amsalpha+}

\begin{abstract}
The main purpose of this paper is to 
make $\overline{C}_{n,n-1}$, 
which is the main theorem of \cite{ka-ma}, 
more accessible. 
\end{abstract}

\maketitle

\section{Introduction}\label{intro}
In spite of its importance, the proof 
of $\overline{C}_{n,n-1}$ is not so easy to access for 
the younger generation, including myself. 
After \cite{ka-ma} was published, 
the birational geometry has drastically developed. 
When Kawamata wrote \cite{ka-ma}, the 
following techniques and results are not 
known nor fully matured. 
\begin{itemize}
\item Kawamata's covering trick, 
\item moduli theory of curves, especially, the 
notion of level structures and the existence of 
tautological families, 
\item various notions of singularities such as 
rational singularities, canonical singularities, 
and so on. 
\end{itemize}
See \cite[\S 2]{ka-d}, 
\cite[Section 5]{ak}, 
\cite[Part II]{ao}, \cite{vo}, \cite{v2}, and \cite{km}. 
In the mid 1990s, de Jong gave us fantastic results:~
\cite{d1} and \cite{d2}. 
The alteration paradigm generated the weak semistable 
reduction theorem \cite{ak}. 
This paper shows how to simplify the proof of the main 
theorem of \cite{ka-ma} by using the weak semistable 
reduction. The proof 
may look much simpler than 
Kawamata's original proof (note that 
we have to read \cite{v1} and \cite{v2} to understand \cite{ka-ma}). 
However, the alteration theorem 
grew out from the deep investigation of the moduli 
space of stable pointed curves (see \cite{d1} and 
\cite{d2}). 
So, don't misunderstand the real value of this 
paper. 
We note that we do not enforce Kawamata's 
arguments. We only recover his main result. 
Of course, this paper is not self-contained. 

The following result is the main theorem of \cite{ka-ma}. 
We call this $\overline{C}_{n,n-1}$ in this paper. 
Here, $n$ means the dimension of $X$. 

\begin{thm}[{\cite[Theorem 1]{ka-ma}}]\label{old}
Let $f:X\longrightarrow Y$ be a dominant morphism 
of algebraic varieties defined over the 
complex number field $\mathbb C$. 
Assume that the general fiber $X_y=f^{-1}(y)$ is an irreducible 
curve. 
Then we have the following 
inequality for logarithmic Kodaira dimensions$:$ 
$$
\overline{\kappa}(X)\geq \overline{\kappa}(Y)
+\overline{\kappa}(X_y). 
$$
\end{thm}

In Section \ref{sec2}, we will give a proof to 
\cite[Theorem 2]{ka-ma}, which is stronger than 
$\overline{C}_{n,n-1}$. See the inequality 
$(\overline{C}'_{n,n-1})$ in the first paragraph 
of the proof below. 

Note that our reference list does not cover all the 
papers treating the related topics. 
We apologize in advance to the colleagues 
whose works were not appropriately mentioned in 
this paper. 
In the proof of the main theorem, we do not refer to the original 
results since they are scattered in various papers. 
Mori collected them nicely in \cite[\S 6, \S7]{mori}. 

\begin{ack} 
I was grateful to the Institute for Advanced Study 
for its hospitality.  
I was partially supported by a grant from the 
National Science Foundation:~DMS-0111298. 
I would like to thank Professor Noboru Nakayama for comments 
and Professor Kalle Karu for giving me \cite{karu}. 
\end{ack}

\begin{notation}
We will work over $\mathbb C$ 
throughout this paper. 
For the basic properties of the logarithmic 
Kodaira dimension, see \cite{ori}, \cite{mont}, 
\cite{hon}, and \cite[\S 1]{ka-ma}. 
\begin{enumerate}
\item[(i)] Let $X$ be a (not necessarily complete) variety. 
Then $\overline{\kappa}(X)$ denotes the {\em{logarithmic 
Kodaira dimension}} of $X$. 
\item[(ii)] Let $f:X\longrightarrow Y$ be a dominant morphism between 
varieties and $D$ a $\mathbb Q$-divisor on $X$. 
We can write $D=D_{\xhor}+D_{\xver}$ such 
that every irreducible component of $D_{\xhor}$ (resp.~
$D_{\xver}$) is mapped (resp.~not mapped) onto 
$Y$. 
If $D=D_{\xhor}$ (resp.$D=D_{\xver}$), $D$ 
is said to be {\em{horizontal}} (resp.~{\em{vertical}}). 
\item[(iii)] Let $f:X\longrightarrow Y$ be a birational morphism. 
Then $\Exc(f)$ denotes the exceptional locus of $f$. 
\end{enumerate}
\end{notation}

\section{$\overline{C}_{n,n-1}$}\label{sec2}

Here, we prove the following theorem. 
It is easy to see that this statement is equivalent to 
Theorem \ref{old} by the basic properties of 
the logarithmic Kodaira dimension. 

\begin{thm}[$\overline{C}_{n,n-1}$]\label{new}
Let $f:X\longrightarrow Y$ be a surjective 
morphism with connected fibers between 
non-singular projective varieties $X$ and $Y$. 
Let $C$ and $D$ be simple normal crossing divisors 
on $X$ and $Y$. 
We put $X_0:=X\setminus C$ and $Y_0:=Y\setminus D$. 
Assume that $f(X_0)\subset Y_0$. 
Then 
$$
\overline{\kappa}(X_0)\geq \overline{\kappa}(Y_0)
+\overline{\kappa}(F_0),  
$$ 
where $F_0$ is a sufficiently general fiber of $f_0:=
f|_{X_0}:X_0\longrightarrow Y_0$. 
\end{thm}

Before we start the proof, let us recall the following trivial 
lemma. We will frequently use it without mentioning it. 

\begin{lem}
Let $X$ be a complete normal variety. 
Let $D_1$ and $D_2$ be $\mathbb Q$-Cartier 
$\mathbb Q$-divisors on 
$X$. 
Assume that $D_1\geq D_2$. 
Then $\kappa (D_1)\geq \kappa (D_2)$. 
\end{lem}

\begin{proof}[Proof of {\em{Theorem \ref{new}}}]  
By \cite{ka-ma}, it is sufficient to 
prove 
\begin{equation}\tag{$\overline{C}'_{n.n-1}$}
\kappa (K_X+C-f^*(K_Y+D))\geq \overline{\kappa}(F_0). 
\end{equation}
\begin{step}
By Theorem 2.1 in \cite{ak} (see also \cite[Chapter 2, 
Remark 4.5 and Section 9]{karu}), 
we have the following commutative diagram: 
$$
\begin{matrix}
X&\longleftarrow &X'&\supset &U_{X'}\\
\downarrow& & \downarrow& & \downarrow\\
Y&\longleftarrow &Y'&\supset &U_{Y'}\\
\end{matrix}
$$ 
such that $p:X'\longrightarrow X$ and $q:Y'\longrightarrow Y$ 
are projective birational morphisms, 
$X'$ is quasi-smooth (in particular, 
$\mathbb Q$-factorial) 
and $Y'$ is non-singular, the inclusion 
on the right are toroidal embeddings, 
and such that 
\begin{enumerate}
\item $f':(U_{X'}\subset X')\longrightarrow 
(U_{Y'}\subset Y')$ is toroidal and equi-dimensional, 
\item Let $C':=(p^*C)_{\red}$ and 
$D':=(q^*D)_{\red}$. Then 
$C'\subset X'\setminus U_{X'}$ and 
$D'\subset Y'\setminus U_{Y'}$. 
\end{enumerate}
Since 
$$\overline {\kappa} (X_0)=\kappa (K_X+C) 
=\kappa (K_{X'}+C')
$$ and 
$$
\overline {\kappa} (Y_0)=\kappa (K_Y+D) 
=\kappa (K_{Y'}+D'), 
$$ 
we can replace $f:X\longrightarrow Y$ with 
$f':X'\longrightarrow Y'$. 
For the simplicity of the notation, we omit the 
superscript $'$. 
So, we can assume that $f:X\longrightarrow Y$ 
is toroidal with the above extra assumptions. 
\end{step}
\begin{step}
By taking a Kawamata's Kummer 
cover $q:Y'\longrightarrow Y$, 
we obtain the following commutative diagram: 
$$
\begin{CD}
X@<{p}<<X'\\
@V{f}VV @VV{f'}V\\
Y@<<{q}<Y'
\end{CD}
$$
such that $f':X'\longrightarrow Y'$ is weakly semistable, 
where $X'$ is the normalization of $X\times _YY'$ 
(see \cite[Section 5]{ak}). 
Note that $X'$ is Gorenstein by \cite[Lemma 6.1]{ak}. 
We put $G:=X\setminus U_X$ and 
$H:=Y\setminus U_Y$. 
Then we have 
$$
K_X+C-f^*(K_Y+D)\geq K_X+C_{\xhor}+G_{\xver}
-f^*(K_Y+H). 
$$
Therefore, we can check that 
$$
p^*(K_X+C-f^*(K_Y+D))\geq K_{X'/Y'}+(p^*C)_{\xhor}. 
$$ 
We note that $(p^*C)_{\xhor}=p^*(C_{\xhor})$. 
So, it is sufficient to prove that 
$\kappa (K_{X'/Y'}+(p^*C)_{\xhor})\geq \overline{\kappa} 
(F_0)$.
\end{step}

\begin{step}
Let $F$ be a general fiber of $f:X\longrightarrow Y$. 
We put $g:=g(F)$:~the genus of $F$. 

\begin{case}[$g\geq 2$] 
In this case, 
$$
\kappa(K_{X'/Y'}+(p^*C)_{\xhor})\geq \kappa(K_{X'/Y'})\geq 1
=\overline{\kappa}(F_0).  
$$ 
The last inequality is well-known. 
See, for example, \cite[(7.5) Theorem]{mori} 
and \cite[Theorem 5.3, Remark 5.4]{can}. 
So, we stop the proof in this case. 
\end{case} 
\begin{case}[$g=1$] 
It is well-known that 
$$
\kappa(K_{X'/Y'})\geq \Var(f')=\Var(f)\geq 0. 
$$ 
See \cite[(7.5) Theorem]{mori} 
and \cite[Theorem 5.3, Remark 5.4]{can}. 
For the definition of the {\em{variation}} $\Var(f)$, 
see, for instance, \cite[p.329]{v3} and \cite[(7.1)]{mori}. 
So, if $C$ is vertical or $\Var(f)\geq 1$, 
then we obtain 
$$
\kappa(K_{X'/Y'}+(p^*C)_{\xhor})\geq\overline{\kappa}
(F_0).  
$$ 
Therefore, we can assume that $\Var(f)=0$ and 
$C$ is not vertical. 
By Kawamata's covering trick, 
we obtain the following 
commutative diagram: 
$$
\begin{CD}
X'@<{\pi}<<X''\\
@V{f'}VV @VV{f''}V\\
Y'@<<{\eta}<Y'', 
\end{CD}
$$where $\eta:Y''\longrightarrow Y'$ is a 
finite cover from a 
non-singular projective variety $Y''$, $f'':X'':=X'\times _{Y'}Y''\longrightarrow Y''$ is weakly semistable, and $f''$ is birationally equivalent to 
$Y''\times E\longrightarrow Y''$. Here, 
$E$ is an elliptic curve. 
Note that, if we need, we can blow-up $Y'$ and replace $X'$ 
with its base change before taking the cover. 
It is because the property of a morphism being weakly semistable 
is preserved by a base change under some mild conditions 
(cf.~\cite[Lemma 6.2]{ak}). 
For details, see \cite[Lemma 6.2]{ak} and the proof of 
\cite[Corollary 19]{ka-d}. 
Since 
$$
\pi^*(K_{X'/Y'}+(p^*C)_{\xhor})=K_{X''/Y''}+\pi^*((p^*C)_{\xhor}), 
$$ 
it is sufficient to prove 
$\kappa (K_{X''/Y''}+\pi^*((p^*C)_{\xhor}))\geq 1$. 
Let $\alpha:\widetilde X\longrightarrow Y''\times E$, 
$\beta:\widetilde X\longrightarrow X''$ be a common 
resolution. 
Since $X''$ has only rational Gorenstein singularities, 
$X''$ has at worst canonical singularities. Thus, 
we obtain 
$$
\kappa (K_{X''/Y''}+\pi^*((p^*C)_{\xhor}))
=\kappa (K_{\widetilde{X}/Y''}+\beta^*\pi^*((p^*C)_{\xhor})). 
$$ 
On the other hand, 
$$
K_{\widetilde{X}/Y''}=K_{\widetilde{X}/Y''\times E}
+K_{Y''\times E/Y''}=:A
$$ 
is an effective $\alpha$-exceptional divisor such that 
$\xSupp A=\Exc(\alpha)$. 
Let $B$ be an irreducible component 
of $\beta^*\pi^*((p^*C)_{\xhor})$ such that 
$B$ is dominant onto $Y''$. 
Then 
$$
m(A+\beta^*\pi^*((p^*C)_{\xhor}))\geq \alpha^*\alpha_*B, 
$$ 
for a sufficiently large integer $m$. 
Therefore, if is sufficient to prove 
$\kappa(Y''\times E, \alpha_*B)\geq 1$. 
It is true by \cite[Corollary 5.4]{mad}. 
Thus, we finish the proof when $g=1$. 
\end{case}
\begin{case}[$g=0$]
As in the above case, we can 
take a finite cover and obtain the 
following commutative diagram: 
$$
\begin{CD}
X'@<{\pi}<<X''\\
@V{f'}VV @VV{f''}V\\
Y'@<<{\eta}<Y'', 
\end{CD}
$$ 
where $f''$ is birationally equivalent to 
$Y''\times \mathbb P^1\longrightarrow Y''$. 
We can further assume that 
all the horizontal components of $\pi^*((p^*C)_{\xhor})$ 
are mapped onto $Y''$ birationally. 
\end{case}
\begin{lem}[cf.~{\cite[Section 7]{can}}] \label{a}
Let $f:V\longrightarrow W$ be a surjective morphism 
between non-singular projective varieties 
with connected fibers. 
Assume that $f$ is birationally equivalent to 
$W\times \mathbb P^1 \longrightarrow W$. 
Let $\{C_k\}$ be a set of distinct 
irreducible divisors such that 
$f:C_k\longrightarrow W$ is birational for every $k$ $(k\leq 3)$. 
Then 
$$\kappa (K_{V/W}+C_1+C_2)\geq 0$$ 
and 
$$
\kappa (K_{V/W}+C_1+C_2+C_3)\geq 1.
$$ 
\end{lem}
\begin{proof}
By modifying $V$ and $W$ birationally (see also \cite
[Lemma 7.8]{can}) and replacing $C_k$ with its strict transform, 
we can assume that there exists a simple normal crossing 
divisor $\Sigma$ on $W$ such that 
$$
\varphi_{ij}:V_0:=f^{-1}(W_0)\simeq W_0\times \mathbb P^1 
$$ 
with $\varphi_{ij}(C_i|_{V_0})=W_0\times\{0\}$ and 
$\varphi_{ij}(C_j|_{V_0})=W_0\times\{\infty\}$ for 
$i\ne j$, 
where $W_0:=W\setminus \Sigma$. 
We can further assume that there exists 
$\psi_{ij}:V\longrightarrow \mathbb P^1$ such that 
$\psi_{ij}|_{V_0}=p_2\circ\varphi_{ij}$, 
where $p_2$ is the second projection $W_0\times \mathbb P^1 
\longrightarrow \mathbb P^1$. 
We also assume that $\cup_kC_k\cup(f^*\Sigma)_{\red}$ 
is a simple normal crossing divisor. 
 we obtain 
\begin{eqnarray*}
\wedge {\psi_{ij}}^*\left( \frac{dz}{z} \right)&
\in&\Hom_{\mathcal O_V}
(f^*(K_W+\Sigma), K_V+C_i+C_j+(f^*\Sigma)_{\red})\\
&\simeq& H^0(V,K_{V/W}+C_i+C_j+
(f^*\Sigma)_{\red}-f^*\Sigma)\\ 
&\subset& 
H^0(V,K_{V/W}+C_i+C_j)
\end{eqnarray*} 
for $i\ne j$, 
where $z$ denotes a suitable inhomogeneous 
coordinate of $\mathbb P^1$ (see 
\cite[Lemma 7.12]{can}). 
Therefore, 
$$\dim_{\mathbb C} 
H^0(V, K_{V/W}+C_1+C_2)\geq 1$$ and 
$$\dim _{\mathbb C} H^0(V, 
K_{V/W}+C_1+C_2+C_3)\geq 2.$$ 
Thus, we obtain the required result. 
\end{proof}
Apply Lemma \ref{a} to $\widetilde X\longrightarrow Y''$, 
where 
$\beta:\widetilde X\longrightarrow X''$ is a resolution of $X''$. 
Then we obtain 
$$
\kappa (K_{\widetilde{X}/Y''}+\beta^*\pi^*((p^*C)_{\xhor}))
\geq \overline{\kappa} (F_0). 
$$
\end{step}
Thus, we complete the proof. 
\end{proof}

\ifx\undefined\bysame
\newcommand{\bysame|{leavemode\hbox to3em{\hrulefill}\,}
\fi

\end{document}